\documentclass[12pt]{article}
\usepackage{amssymb,amscd,amsmath,amsthm}
\usepackage{latexsym,amstext}
\usepackage{latexsym,amstext}
\usepackage{color}
\usepackage{graphicx}
\usepackage{epstopdf}
\usepackage{cite}

\begin{document}

\title{\bf An Extended Picard Method to solve non-linear systems of ODE with applications. }

\author{M. Gadella$^1$, L.P. Lara$^{2,3}$ \\ \\
$^1$ Departamento de F\'{\i}sica Te\'orica, At\'omica y Optica  and IMUVA, \\
Universidad de Va\-lladolid, 47011 Valladolid, Spain\\ 
$^2$ Instituto de F\'isica Rosario, CONICET-UNR, \\ 
Bv. 27 de Febrero, S2000EKF Rosario, Santa Fe,  Argentina.\\
$^3$ Departamento de Sistemas, Universidad del Centro Educativo\\
Latinoamericano, Av. Peegrini 1332, S2000 Rosario, Argentina
}

\maketitle

\begin{abstract}

We provide of a method to integrate first order non-linear systems of differential equations with variable coefficients. It determines approximate solutions given initial or boundary conditions or even for Sturm-Liouville problems. This method  is a mixture between an iterative process, a la Picard, plus a segmentary integration, which gives explicit approximate solutions in terms of trigonometric functions and polynomials. The segmentary part is particularly important if the integration interval is large. This procedure provide a new tool so as to obtain approximate solutions of systems of interest in the analysis of chemical reactions. We test the method on some classical equations like Mathieu, Duffing quintic equation or Bratu's equation and have applied it on some models of chemical reactions.

\end{abstract}

\section{Introduction}

For non-linear ordinary differential equations, non-linear ordinary systems or linear equations with variable coefficients, there exists a panoplia of methods for approximate solutions in the case of an exact solution is unknown. Some of the most celebrated methods are Runge-Kutta, Taylor, Harmonic Balance or Picard (that we here denote as Standard Picard for reasons which wi become evident soon) \cite{HS,MIC,CHI,CL,NV,GS,NOR,AF}. For a given problem for which the solution passes for solving a given differential equation, we should use the method that gives more precision using less steps and less CPU times and this depends on the equation. The more resolution methods available to the scientists, the more possibilities for finding better approximate solutions.

In the present paper, we propose a method to find approximate solutions for certain first order non-linear systems of ordinary differential equations that may serve as solutions of certain chemical processes.  In contrast with some one-step methods previously considered by the authors \cite{LL,LG,GL,GL1}, we here construct iterative approximate solutions, an idea inspired in the Picard method (aka Standard Picard method). In fact, this is a generalization and improvement of the Standard Picard method that we believe may be efficient for obtaining approximate local solutions for systems of first order (at least) ordinary differential equation with potential interest in chemical reactions. From now on, we call this generalization of the Picard method as the {\it Extended Picard method}. 

The basic ideas leading to this method are not new and similar ideas have been introduced some time ago, such as those targeted for their applications on second order differential equations \cite{RAM}. Here, we intend to offer the readers a we written presentation of this method accesible to all with solid mathematical background. At the same time we include some improvements of previous contributions on the subject such as a proof of the convergence of the iterations to the exact solution that somehow generalizes the proof given in \cite{RAM} for the case thereby studied. In addition, we add some applications including one that may have interest in the study of chemical reactions.

In some practical situations, particularly when the integration interval is too large, Extended Picard, as any other iterative procedure, may require many iterations for a reasonable precision far away from the point determining the initial value (in the case of first order systems). However, this enhances the possibility of arithmetic errors or other deviations. This is why that for these situations, we need to complete the procedure with a segmentary integration that permits a much smaller number of iterations. In our opinion this is another improvement that permits the implementation of the method using class-room computers. This paper is a kind of version of a unpublished preprint \cite{PRE}.

This paper is organized as follows: An introduction to the Extended Picard method is given in Section 2, where we also show the uniform convergence of the sequence of iterations to the exact solution. We also explain how to complete the iteration method with segmentary integration. Testing the method on some classical differential equations with wide applications is the objective of Section 3. In Section 4, we make an extensive application of Extended Picard to two equations describing chemical reactions: Glycolisis and Brusselator. We compare the precision of our method with some others like Taylor or Standard Picard. In the case of Brusselator, we show how a good choice of the form of the equation, improves properties of the solution, such as precision. We close the paper with some Concluding remarks.

\section{Position of the problem and convergence of the approximate solutions.}

Let us consider a $N$-dimensional system of first order non-linear differential equations, written in this form:

\begin{equation}\label{1}
\mathbf y'(x) =\mathbf F(x,\mathbf y)\,, \qquad \mathbf y \equiv \left(\begin{array}{c} y^1(x) \\ y^2(x) \\\vdots \\y^N(x) \end{array} \right)\,, \qquad F(x,\mathbf y) \equiv \left( \begin{array}{c} f^1(x,\mathbf y) \\ f^2(x,\mathbf y) \\ \vdots \\ f^N(x,\mathbf y) \end{array} \right) \,, 
\end{equation}
where $F(x,\mathbf y)$ is continuous as a function $\mathbf F(x,\mathbf y): [a,b] \times \mathbb R^N \longmapsto \mathbb R^N$ and satisfies a Lipschitz condition with respect to the last $N$ variables. Here $[a,b]$ is a real interval. Equation \eqref{1} may always be written as

\begin{equation}\label{2}
\mathbf y'(x) = A\cdot \mathbf y(x) + \mathbf G(x,\mathbf y)\,, \qquad  \mathbf G(x,\mathbf y) := \mathbf F(x,\mathbf y)-A \cdot \mathbf y\,,
\end{equation}
where $A$ is a constant $N\times N$ matrix. Obviously, $\mathbf G(x,\mathbf y)$ satisfies the Lipschitz condition with respect to the $\mathbf y$ variables. 

Let $\mathbf y_0(x)$ be a particular solution of the homogeneous system $\mathbf y'(x)=A\cdot \mathbf y(x)$ satisfying certain initial or boundary conditions on the interval $[a,b]$. Now, let $\mathbf y_k(x)$ be the solution of the system

\begin{equation}\label{3}
\mathbf y'_k(x)= A \cdot \mathbf y_k(x) +\mathbf G(x,\mathbf y_{k-1}(x))\,,\qquad k=1,2,\dots\,,
\end{equation}
under the same boundary conditions we have imposed to $\mathbf y_0(x)$.  The explicit form of the solution $\mathbf y_k(x)$ is

\begin{equation}\label{4}
\mathbf y_k(x) =\mathbf y_0(x) + \mathbf d_{k-1}(x)\,,
\end{equation}
where $\mathbf d_{k-1}(x)$ is the particular solution of the previous iteration, which can be written as

\begin{equation}\label{5}
\mathbf d_{k-1}(x) = \int_a^x e^{(x-s)A}\,\mathbf G(s,\mathbf y_{k-1}(s))\,ds\,.
\end{equation}

Once we have obtained the iterations in \eqref{4}, we need to show that the sequence $\{\mathbf y_k(x)\}$ converges uniformly to a solution of \eqref{1} on the interval $[a,b]$. 

\subsection{On the issue of convergence}

Prior to this discussion, we need to recall some simple facts concerning the norm of a matrix. Although all norms on a space of matrices give the same topology and, hence, the same notion of convergence, some norms are more convenient than other, depending in our purpose. Assume that $B$ is a $N\times N$ square matrix and that $\mathbf x$ is an arbitrary vector in $\mathbb R^N$. Then, we may define the norm, $||B||$, of $B$ as

\begin{equation}\label{6}
||B||:= \sup_{||\mathbf x||=1}||B\mathbf x||\,,
\end{equation}
where to any $\mathbf x =(x_1,x_2,\dots,x_N) \in \mathbb R^N$, we assign its usual norm, $||\mathbf x||:= \sqrt{x_1^2 +x_2^2+\dots +x_N^2}$. 
Norm \eqref{6} has an interesting property, which is \cite{BN}

\begin{equation}\label{7}
||B^n|| \le ||B||^n\,,\qquad n=1,2,\dots\,.
\end{equation}

Now, assume that $|t|\le c$. Then,

\begin{eqnarray}\label{8}
||e^{tA}|| =\left|\left| \sum_{n=0}^\infty \frac{t^n\,A^n}{n!}   \right|\right| \le \sum_{n=0}^\infty \frac{|t|^n \,||A||^n}{n!} \le \sum_{n=0}^\infty \frac{c^n\,||A||^n}{n!} = e^{c\,||A||} \,,
\end{eqnarray}
where in the first inequality in \eqref{8}, we have used \eqref{7} plus the continuity of the norm.

Let us consider the iterations of the form \eqref{4}, where we first choose $k=1$. We estimate an upper bound for the norm in $\mathbb R^N$, valid for any $x\in [a,b]$ and $c:= \max\{|a|,|b|\}$:

\begin{equation}\label{9}
|| \mathbf y_1(x) - \mathbf y_0(x)|| \le \int_a^x || e^{(s-x)A}\, \mathbf G(s,\mathbf y_0(s)||\,ds \le e^{2c\,||A||} \int_a^x ||\mathbf G(s,\mathbf y_0(s)||\,ds\,.
\end{equation}

The continuity of the functions $\mathbf G(x,\mathbf y)$ and of $\mathbf y(x)$ on the interval $[a,b]$, shows that the norm under the last integral sign in \eqref{9} must be bounded by a certain positive constant, say $H$. We denote by $M:= e^{2c\,||A||}$ so that the inequalities in \eqref{9} give

\begin{equation}\label{10}
|| \mathbf y_1(x) - \mathbf y_0(x)|| \le MH\,(x-a)\,.
\end{equation}

The second step is to find an upper bound for

\begin{eqnarray}\label{11}
|| \mathbf y_2(x) - \mathbf y_1(x)|| \le \int_a^x || e^{(s-x)A} \,  [\mathbf G(s,\mathbf y_{1}(s)) - \mathbf G(s,\mathbf y_{0}(s))]\,||\,ds \nonumber \\ [2ex] \le M \int_a^x || \mathbf G(s,\mathbf y_{1}(s)) - \mathbf G(s,\mathbf y_{0}(s))||\,ds\,.
\end{eqnarray}

By hypothesis, the function $\mathbf G(x,\mathbf y)$ satisfies the Lipschitz condition with respect to the last $N$ variables. This means that, there exists a positive constant $K>0$, such that

\begin{equation}\label{12}
||\mathbf G(s,\mathbf y_{k-1}(x)) - \mathbf G(s,\mathbf y_{k-2}(x))|| \le K \,||\mathbf y_{k-1}(x)-\mathbf y_{k-2}(x)||\,, \qquad k=2,3,\dots\,.
\end{equation}

Let us use  \eqref{12}, with $k=2$, in the last integral in \eqref{11}. This gives:

\begin{eqnarray}\label{13}
|| \mathbf y_2(x) - \mathbf y_1(x)|| \le MK \int_a^x || \mathbf y_1(s)- \mathbf y_0(s)||\,ds \le HM^2K^2 \,\frac{(x-a)^2}{2}\,.
\end{eqnarray}

By induction, one finds

\begin{equation}\label{14}
||\mathbf y_k(x)- \mathbf y_{k-1}(x)|| \le HM^k K^k\,\frac{(x-a)^k}{k!} \le HM^k K^k\,\frac{(b-a)^k}{k!}  \longmapsto 0\,, 
\end{equation} 
if $k\longmapsto\infty$, for any value $x\in [a,b]$. This implies that the sequence \eqref{4} {\it converges uniformly to some continuous function} on the interval $[a,b]$. Let us call $\mathbf y(x)$ to such a limit.

The first derivative of the functions in the sequence $\{\mathbf y_k(x)\}$, in \eqref{4}, obviously exists. If we apply the Leibnitz formula to \eqref{5}, we have

\begin{equation}\label{15}
d'_k(x) = G(x,\mathbf y_{k-1}(x)) + \int_a^x e^{(x-a)A}\,A \, G(s,\mathbf y_{k-1}(s))\,ds\,,
\end{equation}
where the prime denotes derivative with respect to $x$.

In order to show the uniform convergence of the sequence of the derivatives, $\{ \mathbf y_k(x)\}$, we use a similar argument as above and write, 

\begin{eqnarray}\label{16}
\mathbf y_{k+1}'(x) - \mathbf y_k'(x) = \mathbf G(x, \mathbf y_k(x)) - \mathbf G(x, \mathbf y_{k-1}(x)) \nonumber \\[2ex] + \int_a^x e^{(x-s)\,A}\, A \, [\mathbf G(s, \mathbf y_k(s))- \mathbf G(s, \mathbf y_{k-1}(s))]\,ds\,.
\end{eqnarray}

Then, taken norms, using the Lipschitz condition for $\mathbf G(x, \mathbf y)$, we have that 

\begin{eqnarray}\label{17}
||\mathbf y_{k+1}'(x) - \mathbf y_k'(x)|| \le K\,|| \mathbf y_k(x) - \mathbf y_{k-1}(x)|| + KM \,||A|| \int_a^x || \mathbf y_k(s) - \mathbf y_{k-1}(s)||\,ds  \nonumber \\ [2ex] \le HM^k K^{k+1} \, \frac{(b-a)^k}{k!} + H\,||A||\, M^{k+1}\, K^{k+1}\, \frac{(b-a)^{k+1}}{k!} \longmapsto 0\,, \qquad {\rm if} \;\; k \longmapsto \infty\,,
\end{eqnarray}
where we have used \eqref{14} in the second inequality in \eqref{17}. This shows the uniform convergence of the sequence of the derivatives. 

This shows that the limit function $\mathbf y(x)$ is differentiable and

\begin{equation}\label{18}
\mathbf y'(x) = \lim_{n\to\infty} \mathbf y'_n(x)
\end{equation}
on $[a,b]$. 

Then, let us go back to \eqref{5} and take the limits as $k\longmapsto \infty$. Using the Lebesgue dominated convergence theorem and the continuity of $\mathbf G(x, \mathbf y)$, we obtain after the limit,

\begin{equation}\label{19}
\mathbf y(x) = \mathbf y_0(x) + \int_a^x e^{(x-s)\,A}\, \mathbf G(s,\mathbf y(s))\,ds\,.
\end{equation}

Then, if we take derivatives in \eqref{19} following the Leibnitz rule, we have

\begin{eqnarray}\label{20}
\mathbf y'(x) = \mathbf y'_0(x) + \mathbf G(x, \mathbf y(x)) +A \int_a^x e^{(x-s)\,A} \, \mathbf G(s,\mathbf y(s))\,ds  \nonumber \\[2ex] = A \left[ \mathbf y_0(x) + \int_a^x e^{(x-s)\,A} \, \mathbf G(s,\mathbf y(s))\,ds \right] + \mathbf G(x, \mathbf y(x)) = A \mathbf y(x) + \mathbf G(x, \mathbf y(x))\,,
\end{eqnarray} 
so that $\mathbf y(x)$ satisfies \eqref{2}. Moreover, the solution $\mathbf y(x)$ as well as all the functions in the sequence $\{\mathbf y_k(x)\}$ satisfy the same initial condition at $x=a$ than $\mathbf y_0(x)$. 

\subsection{Extension of approximate solutions: Segmentary integration.}

Once we have obtained an approximate solution, we may use it for practical purposes. However, this procedure has its own limitations. First, the bigger the interval is, the bigger number of iterations we need in order to get an acceptable approximation. The bigger the number of iterations, the bigger the number of various errors. These  errors play a role similar to a noise, the precision that one may reach by increasing the number of iterations is blurred by these numerical errors. In addition, the iterative integration should be performed by a software of symbolic calculus.

In order to overcome these inconveniences, we may resort to a combination of iterative integration and segmentary integration or one step integration. It goes as follows: Let us divide the interval of integration on a certain number of subintervals of equal length (this is not necessary, although we choose equal length here for practical convenience). The points the separate these subintervals or {\it nodes} give a sequence $x_0=a, x_1,\dots,x_k, \dots,x_m=b$. Now, Let us assume that $\mathbf y(x)$ is the exact solution, defined as above and obtained by the previous procedure on the interval $[x_0,x_1]$. Its value at $x_1$ is $\mathbf y(x_1)$. However, we do not know its precise value in general, so that it may be invalid as initial condition for the solution on the second interval $[x_1,x_2]$. Then, we proceed as follows: We start with a seed solution verifying the initial condition to start with the process. Then, assume that $\mathbf y_p(x)$ is the approximate solution, obtained through the iteration process, that we have chosen on the first interval $[x_0,x_1]$. Then, we determine the initial value on $[x_1,x_2]$ at the point $x_1$ to be equal to $\mathbf  y_p(x_1)$. Observe that, due to the uniform convergence of the approximate solutions to the exact solution, the error consisting in taking $\mathbf  y_p(x_1)$ as initial value for the second interval instead of $\mathbf y(x_1)$ is controlled.  Then, we repeat all the above procedure on the interval  $[x_1,x_2]$.

It should be clear that the exact solution on $[x_1,x_2]$ does not match with the exact solution on $[x_0,x_1]$, due to the fact that we have not been able to use $\mathbf y(x_1)$ as the initial condition for $[x_1,x_2]$. The seed solution for $[x_1,x_2]$ is the solution obtained on $[x_0,x_1]$ with the desired precision. We do this for all segments, in the sense that the seed solution for $[x_{s-1},x_s]$ is the final solution chosen for $[x_{s-2},s_{s-1}]$.     Nevertheless, the approximate solution is always continuous on $x_1$ by construction, no matter which interaction we have chosen for the approximate solution on $[x_1,x_2]$. 

We repeat this procedure on each of the intervals $[x_{s-1},x_s]$ so as to obtain approximate solutions. As we go further in $s$, difficulties to obtain the solution arose, as mentioned before. To overcome these difficulties, let us go back to \eqref{3} and approximate the non-linear terms in $\mathbf G(x,\mathbf y_{k-1}(x))$ by minimal squares with polynomials of order one or three. This latter choice depends on which one gives better precision, something that should be checked by numerical experiments and that may depend on the explicit form of $\mathbf G(x,\mathbf y_{k-1}(x))$.

\section{Testing the method}

We test the method on three classical equations with a wide range of applications: Mathieu equation, quintic Duffing equation and Bratu equation.  

\subsection{The Mathieu equation}

The Mathieu equation is the simplest differential equation of Hill type and is customary to write it on the following form:

\begin{equation}\label{21}
y''(x) +(r-2q\,\cos (2x))\,y(x)=0\,,
\end{equation}
where $r$ is often called the characteristic value of eigenvalue, which is usually determined through some boundary conditions and $q$ is a fixed data. In general, one looks for periodic solutions with boundary conditions at some finite points, say 0 and $p$, where $p$ is the period, of the form $y(0)=y(p)$, $y'(0)=y'(p)$, usually $p=\pi$ or $p=2\pi$. The literature concerning the properties of the Mathieu equation and, its solutions, the Mathieu functions is wide. Let us just cite  \cite{MAC,FAR,BC,ZIE,SHI,GGL}. Typically, each particular solution of \eqref{21} could be written in terms of  one of the Mathieu functions $C(r,q,x)$ and $S(r,q,x)$, where $r$ is one of the possible eigenvalues, being given the value of $q$. See \cite{ABR}.

Equation \eqref{21} may be written in terms of a first order two dimensional system in the form,

\begin{eqnarray}
y'(x) & = & z(x)\,, \label{22} \\[2ex]  z'(x) & = & -(r-2q\,\cos(2x))\,y(x)\,. \label{23}
\end{eqnarray}

Let us impose the initial conditions, $y(0)=1$ and $z(0)=0$ and choose $q=0.05$.  One  particular solution is 

\begin{equation}\label{24}
y_M(x) = (0.81413- 0. 571483\,i) \, C(1,0.05,x)\,,
\end{equation}
where the subindex $M$ stands for Mathieu. 

Let us use the extended Picard method to solve the system \eqref{22}-\eqref{23} using the seeds functions

\begin{equation}\label{25}
y_0(x)= \cos x\,, \qquad z_0(x) = -\sin x\,. 
\end{equation}

After \eqref{3}, we obtain the following recurrence relations for $k=1,2,\dots$,

\begin{eqnarray}
y_k(x) \cos x + r \int_0^x (\sin (t+x) -\sin(3t-x))\, y_{k-1}(t)\,dt\,, \label{26} \\[2ex]  z_k(x) =-\sin x +r \int_0^x (\cos(t+x) + \cos (3t-x))\, z_{k-1}(t)\,dt\,. \label{27}
\end{eqnarray}

Since we have already an exact solution given by \eqref{24}, we may compare the precision of the $n$-th approximation to the exact solution. As usual, we define the error as

\begin{equation}\label{28}
e:= \frac 1{2\pi} \int_0^{2\pi} (y_M(x)-y_n(x))^2\,dx\,.
\end{equation}

We may appreciate the precision of the method just with one iteration. For $n=1$, the error \eqref{28} gives $e=3.33 \,10^{-5}$, while for $n=2$ is $e=3.57\, 10^{-8}$. For $n>2$, we obtain $e< 10^{-20}$. We also may write the approximate solution of \eqref{26}-\eqref{27} with just two iterations, which is

\begin{eqnarray}\label{29}
y_2(x) = (1.00639 + 0.0003125\, x^2)\,\cos x - 0.00640625 \, cos 3x \nonumber\\[2ex]  + 0.0000130208\, \cos 5x + 0.0246875 \, x\sin x- 0.00015625 \, x\sin 3x\,.
\end{eqnarray}

\begin{eqnarray}\label{30}
z_2(x) = 0.0253125\, x\cos x -0.00046875\,x\cos 3x - 0.981706 \, \sin x \nonumber\\[2ex] - 0.0003125\, x^2\sin x + 0.0190625\, \sin 3x - 0.0000651042\, \sin 5x\,.
\end{eqnarray}

It is noteworthy that already two iterations give a good approximation to the exact solution (see Figure 1), which cannot be expressed in terms of elementary functions. 

\begin{figure}
\centering
\includegraphics[width=0.5\textwidth]{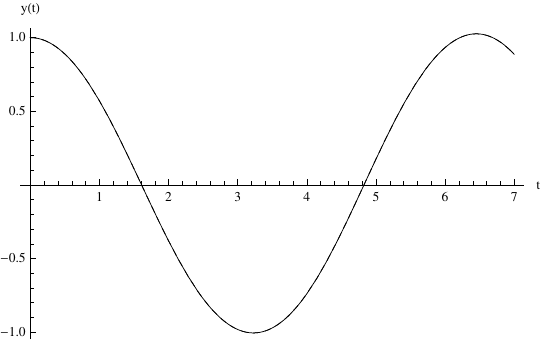}
\caption{\small Graphics for \eqref{29}. This is the approximate iterative solution for $n=2$. It is noteworthy how good the exact solution matches with ours.
\label{Figure1}}
\end{figure}

Our aim next is to compare the results due to the Extended Picard Method to the segmentary Taylor method at different orders. Let us underline that the aim is not to compare the segmentary method with more sophisticated multistep method such as \cite{AF}, but to compare the simplest one step methods with ours. 

 As is well known, the local error produced by the Taylor method is of order of $o(h^n)$, which is incidentally the same error produced by Runge-Kuta, where $h$ is the segment length. We integrate on the interval $[0,2\pi]$. The errors for two different segment lengths are displayed in the next table, where $T_k$ means Taylor method of order $k$.

\vskip1cm

\centerline{$
\begin{array}
[c]{ccccc}
h &T_2 & T_3 & T_4  & T_5  \\[2ex]
0.1 & 1.48\,10^{-5} & 1.62\, 10^{-8} & 4.20\,10^{-12} & 6.57 \,10^{-15}  \\[2ex]
0.5 & 9.93 \, 10^{-3} & 2.62\,10^{-4} & 1.50\, 10^{-6} & 4.70\, 10^{-8}
\end{array}
$}

\bigskip

TABLE 1.- Relative errors produced by the integration of \eqref{21} with the given conditions using the Taylor method with different orders. 

\vskip1cm

Next, let us add the segmentary integration as introduced in Section 2.2. Let us use polynomials of degree one or three on each of the segments and two different values, $h=0.1$ and $h=0.5$ for the segment length. The results are displayed in the following table:

\vskip1cm

\centerline{$
\begin{array}
[c]{ccccccccc}
h &{\rm iteration} & 1 & 3   & | & h & {\rm iteration} & 1 & 3  \\[2ex]
0.1 & 2 & 6.00\, 10^{-11} & 1.30\,10^{-12} & | & 0.5 & 2 & 8.02\,10^{-9} & 7.24\,10^{-11}  \\[2ex]
0.1 & 3 & 6.00\, 10^{-11}  & 1.27\,10^{-12} & | & 0.5 & 3 &  8.00\,10^{-9} & 7.53\, 10^{-11} \\[2ex]
0.1 & 4 & 6.00\, 10^{-11}  & 1.27\,10^{-12} & | & 0.5 & 4 & 8.00\,10^{-9} & 7.53\, 10^{-11} \\[2ex]
0.1 & 5 &  6.00\, 10^{-11}  & 1.27\,10^{-12} & | & 0.5 & 5 & 8.00\,10^{-9} & 7.53\, 10^{-11}
\end{array}
$}

\bigskip

TABLE 2.- Errors obtained by using Extended Picard with segmentary integration using polynomials of degree 1 and 3. We observe that, up to three significative digits, the errors are practically invariant with the number of iterations up to five iterations. 

Now, let us perform another discussion on equation \eqref{21}. Let us determine using Extended Picard the approximate periodic solutions with period $p=\pi$. With these solutions, being given the value of $q$, we may determine  the possible values of the {\it eigenvalue of characteristic value} $r$. In \cite{ABR}, we find the first values of $r$ in terms of $q$, which are

\begin{eqnarray}\label{31}
r_1 = 1-q -\frac{1}{8}\, q^2  + \frac{1}{64}\, q^3 - \frac{1}{1536}\, q^4 + \dots \,, \nonumber\\[2ex]
r_2= 4 + \frac{5}{12}\, q^2 - \frac{763}{13824}\, q^4 + \frac{1002401}{79626240}\, q^6 + \dots \,, \nonumber\\[2ex]
r_3= 9 +\frac{1}{16}\, q^2 - \frac{1}{64}\, q^3 + \frac{13}{20480}\, q^4 + \dots\,,  \nonumber\\[2ex]
r_4= 16 + \frac{1}{30}\, q^2 + \frac{433}{864000}\, q^4 - \frac{5701}{2721600000}\, q^6 + \dots\,,  \nonumber\\[2ex]
r_5= 25 + \frac{1}{48}\, q^2 + \frac{11}{774144}\, q^4 - \frac{1}{147456}\, q^5 + \dots\,.
\end{eqnarray}

Let us apply Extended Picard, for which the iterated equations are ($k=1,2,\dots$)

\begin{equation}\label{32}
y'_k(x) = z_k(x)\,,\qquad z'_k(x) = -r y_k(x) -2q\, \cos(2x)\, y_{k-1}(x)\,,
\end{equation}
with initial conditions $y_k(0)=0$, $z_k(0)=1$. This suggest to use the initial seed given by $y_0(x)\equiv 0$ and $z_0(x) \equiv 1$. We choose the value $q=0.1$ and perform $n$ iterations. The characteristic values are the roots of the equation $y_n(\pi)=0$. In order to estimate the error, we calculate the porcentual relative variation of $r$ obtained after the $n$-th iteration, with $q=0.1$, with respect to the analytic value given by \eqref{31} within the given approximation. For the first three iterations, we obtain the errors listed on Table 3.

\vskip1cm

\centerline{$
\begin{array}
[c]{cccccc}
{\rm iteration} &r_1 & r_2 & r_3  & r_4 & r_5  \\[2ex]
1 & 11 & 1.0\, 10^{-1} &  6.8 \,10^{-3} & 2.0\, 10^{-3} & 8.3\, 10^{-4}  \\[2ex]
2 & 1.4 \, 10^{-1} & 1.0\,10^{-1} & 6.7\, 10^{-3} & 2.0\, 10^{-3} & 8.3\, 10^{-4}  \\[2ex]
3 & 3.4\, 10^{-3} & 1.2\, 10^{-1} & 8.4\, 10^{-4} & 1.4\, 10^{-5} & 3.5\, 10^{-5}
\end{array}
$}

\bigskip

TABLE 3.- Error obtained when calculating the first five characteristic values for one, two or three iterations with respect to the values obtained using relations \eqref{31}. Here, we have used $q=0.1$. 

\vskip1cm

We see that Extended Picard shows a good precision even for a low number of iterations.

\subsection{Quintic Duffing equation}

Non linear oscillators are usually modelled using Duffing type equations \cite{NAG,DAV}. In this Subsection, we apply Extended Picard to the quintic Duffing equation, which written in the form of a two dimensional first order differential system has the form,

\begin{equation}\label{33}
y'(x)= z(x)\,, \qquad z'(x)= -y(x)-a [y(x)]^5\,,
\end{equation}
with $a>0$. We assume that the variable $x$ runs out a finite interval, say $[0,2\pi]$. In order to find a particular solution, we need to fix initial conditions and we choose the very simplest $y(0)=1$ and $y'(0)=0$. For the constant $a$, let us fix $a=1/2$. Then, we integrate with the Extended Picard using the seed

\begin{equation}\label{34}
y_0(x)= \cos x\,, \qquad z_0(x)= -\sin x\,.
\end{equation}

This gives for each $k=1,2,\dots$, the following iterate approximate solutions:

\begin{eqnarray}
y_k(x) & = & \cos x +a \int_0^x [y_{k-1}(t)]^5\,\sin(t-x)\,dt\,, \label{35} \\[2ex]
z_k(x) & = & -\sin x -a \int_0^x [y_{k-1}(t)]^5\, \cos(t-x)\,dt\,. \label{36}
\end{eqnarray}

Obviously, $y'_k(x)=z_k(x)$ for $k=1,2,\dots$. Nevertheless, in order to implement the recursive integration, it is desirable the use of (\ref{35}-\ref{36}) and avoid the use of this derivative. In Figure 2, we compare the solutions obtained with Extended Picard with $k=2$ and eight order Runge-Kutta. Near the initial point, in this case $x=0$, both solutions match fairly well, but at some distance they clearly split and the Extended Picard solution finally blows up. The reason for this behaviour of the Extended Picard solution lies in the fact that this method gives good accurate solutions locally, only. That is to say that Extended Picard is an excellent approximation near the initial conditions. In order to obtain a better accuracy, we have to increase the number of iterations, which increases CPU times and, furthermore, the appearance of numerical errors. 

This is when the use of the segmentary version of the Extended Picard method plays a role in order to improve the approximate solution. Since the precision falls down far away of the initial point, let us make the comparison with the Taylor method at various orders using a larger interval, say $[0,7]$. As the reference solution, we use the eight order Runge-Kutta solution. We define this error as customary:

\begin{equation}\label{37}
e:= \frac 17 \int_0^7 (y_{RK}(x)-y_n(x))^2\,dx\,,
\end{equation}
where $y_{RK}(x)$ is the eight order Runge-Kutta solution and $y_n(x)$ is either the segmentary Extended Picard or the Taylor of $n$-th order. Again, $h$ denotes the interval length. In Table 4, we give the errors produced by the use of Taylor method and in Table 5, the errors obtained after the  segmentary Extended Picard.

\begin{figure}
\centering
\includegraphics[width=0.5\textwidth]{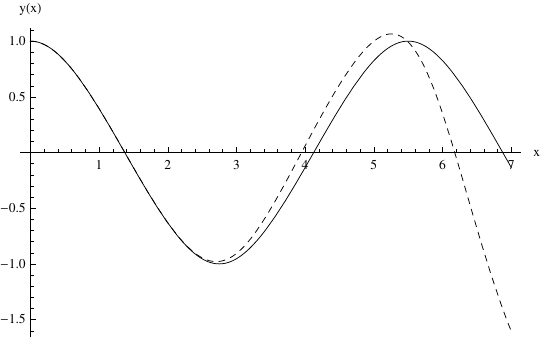}
\caption{\small Comparison of the solution obtained by Extended Picard with $k=2$ (dashed curve) with the solution obatained by eight order Runge-Kutta (continuous curve).
\label{Figure2}}
\end{figure}

\medskip

\centerline{$
\begin{array}
[c]{ccccc}
h &T_2 & T_3 & T_4  & T_5  \\[2ex]
0.1 & 2.55\,10^{-4} & 1.46\, 10^{-6} & 5.170\,10^{-11} & 4.93 \,10^{-12}  \\[2ex]
0.5 & 1.14 & 1.93\,10^{-2} & 1.13\, 10^{-3} & 2.20\, 10^{-5}
\end{array}
$}

\bigskip

TABLE 4.- Relative errors obtained when using the Taylor method of various orders in comparison to the eight order Runge-Kutta solution.

\vskip1cm

\centerline{$
\begin{array}
[c]{ccccccccc}
h &{\rm iteration} & 1 & 3   & | & h & {\rm iteration} & 1 & 3  \\[2ex]
0.1 & 2 & 8.39\, 10^{-9} & 2.19\,10^{-11} & | & 0.5 & 2 & 8.76\,10^{-6} & 1.20\,10^{-6}  \\[2ex]
0.1 & 3 & 8.28\, 10^{-9}  & 1.82\,10^{-11} & | & 0.5 & 3 &  4.45\,10^{-6} & 7.90\, 10^{-8} \\[2ex]
0.1 & 4 & 8.28\, 10^{-9}  & 1.82\,10^{-11} & | & 0.5 & 4 & 4.44\,10^{-6} & 7.91\, 10^{-8} \\[2ex]
0.1 & 5 &  8.28\, 10^{-9}  & 1.82\,10^{-11} & | & 0.5 & 5 & 4.44\,10^{-6} & 7.90\, 10^{-8}
\end{array}
$}

\bigskip

TABLE 5.- Same errors as shown on Table 2 referred to \eqref{31}. 

\vskip1cm

\begin{figure}
\centering
\includegraphics[width=0.5\textwidth]{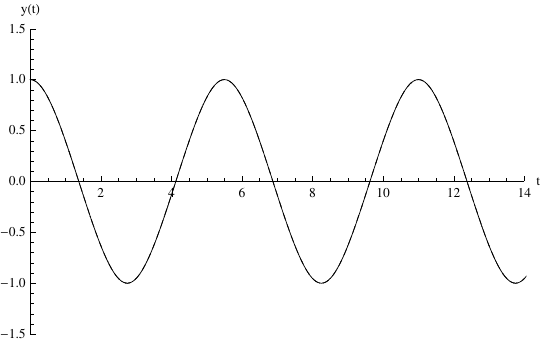}
\caption{\small  Approximated solution of the Duffing equation with segmentary Extended Picard. The Runge-Kutta solution matches with our solution. The interval has been extended to $[0,14]$ in order to display the  periodicity. This periodicity shows since $a>0$.
\label{Figure3}}
\end{figure}

We have shown that the Extended Picard reinforced with segmentary integration gives an excellent approximation to solutions of a strongly non-linear ordinary system such as \eqref{33}. Same can be said for the Mathieu equation as an inspection of Table 2. 

In Figure 3, we have shown one approximate solution of the quintic Duffing equation over two intervals of length equal to the period 7. The solution looks periodic if we restrict ourselves to two periods. This shows the strength of combining Extended Picard with segmentary integration, which is an important conclusion of the present work. 

\subsection{The Bratu equation}

The Bratu equation has arisen in several problems such as thermal reaction processes, a model of the expansion of the Universe, chemical reaction theory, radiative heat transfer and nano-technology \cite{LLI,WAZ,BUC,MD,HE,CAG}. The Bratu equation has the form

\begin{equation}\label{38}
y''(x)+\alpha \,e^y=0\,,
\end{equation}
with $\alpha>0$. We consider the solution with boundary conditions $y(0)=y(1)=0$. In \cite{HS}, we found an analysis towards the determination of an approximated analytical solution to this problem, which uses a method called variational iterational method (VIM). Nevertheless, the exact solution is known to be 

\begin{equation}\label{39}
w(x)= -2 \log \left( \frac{\cosh(0.5(x-0.5)\theta)}{\cosh(0.25\,\theta)} \right)\,,
\end{equation}
where $\theta$ satisfies the following implicit equation \cite{HS}:

\begin{equation}\label{40}
\theta= \sqrt{2\alpha} \, \cosh (0.25\,\theta)\,.
\end{equation}

This problem may have either none, one or two solutions, which are determined  by the value of $\alpha$ with respect to a critical value, $\alpha_c$. Thus, if $\alpha>\alpha_c$, there is no solution, there is one solution if $\alpha=\alpha_c$ and two if $\alpha<\alpha_c$.  This critical value $\alpha_c$ satisfies the following relation:

\begin{equation}\label{41}
4= \sqrt{2 \alpha_c}\, \sinh(0.25\,\theta)\,.
\end{equation}

In \eqref{38}, we may expand $e^y$ on a neighborhood of the origin and leave the quadratic terms,

\begin{equation}\label{42}
y''(x)= -\alpha \left( 1+y + \frac 12\, y^2 \right)\,.
\end{equation}

The variational iterative method (VIM) introduced in \cite{HS} determines the second iteration $y_2(x)$, satisfying $y_2(0)=0$. In order to compare VIM with Extended Picard, we need to give fixed values to constants. If we choose $\alpha=1$, then \eqref{38} gives $\theta= 1.51716459905$. This gives for the second iteration the following expression:

\begin{eqnarray}\label{43}
H_2(x)= kx -\frac{x^2}{2!} - \frac{x^3}{3!} - \frac{(k^2-1) x^4}{4!} + \frac{4x^5}{5!} + \frac{(5k^2-3)x^6}{6!} + \frac{5k(k^2-2) x^7}{7!} \nonumber\\[2ex]  - \frac{25 k^2 x^8}{8!} - \frac{35 k^3 x^9}{9!} - \frac{35 k^4 x^{10}}{10!}\,,
\end{eqnarray}
where $k$ is an arbitrary constant to be determined by a second initial or boundary condition. We choose $k$ such that $y_2(1)=0$ and this gives $k= 0.546936690480377$. We denote the iterative approximate solutions by VIM as $H_n(x)$ in order to distinguish them of the iterations obtained using Extended Picard. 

Let us rewrite  \eqref{42} as a first order two dimensional system as

\begin{equation}\label{44}
y'(x)= z(x)\,, \qquad z'(x)= -\alpha \left( 1+y + \frac 12\, y^2 \right)\,.
\end{equation}

We intend to find an approximate solution of system \eqref{44} using Extended Picard and the initial values $y(0)=0$ and $z(0)=u$, where $u$ is to be determined through the boundary condition $y(0)=1$. As  seed, we use $y_0(x)=0$ and $z_0(x)=u$. The recursive system is obviously, ($k=1,2,\dots,n$)

\begin{equation}\label{45}
y'_k(x) = z_k(x)\,,\qquad z'_k(x) =-1 -y_k(x) - y^2_{k-1}(x)\,,
\end{equation}
which after integration gives

\begin{eqnarray}
y_k(x) =u\, \sin x + \frac 12 \int_0^x (2+ y^2_{k-1}(t))\,\sin(t-x)\,dt\,, \label{46} \\[2ex]
z_k(x) = u\, \cos x - \frac 12 \int_0^x (2+ y^2_{k-1}(t))\,\cos(t-x)\,dt\,. \label{47}
\end{eqnarray}

For instance, take $k=2$ and choose $y_2(1)=0$. then, $u=0.549249$ and

\begin{eqnarray}\label{48}
y_2(x)= - 1.82542 + (1.76722 - 0.274624\, x) \cos x + 0.0581938\, \cos 2x \nonumber\\[2ex] 
+ \left(0.64079 + \frac 12 \, x \right) \, \sin x + 0.0915415 \, \sin 2x\,.
\end{eqnarray}

Since we have an exact solution in \eqref{37}, we can compare the errors obtained using Extended Picard and other methods. For extended Picard and VIM, we shall use the second iteration, for which the solutions are given in \eqref{48} and \eqref{43}, respectively. The other two are the Runge-Kutta solution and Taylor of order 10. As customary, we define the error as 

\begin{equation}\label{49}
e:= \int_0^1 (w(x)- \widetilde y(x))^2\, dx\,,
\end{equation}
where $w(x)$ is the exact solution \eqref{39} and $\widetilde y(x)$ one of the above mentioned approximate solutions. The result is displayed in the next table.

\vskip1cm 

\centerline{$
\begin{array}
[c]{cccc}
{\rm EP} & {\rm VIM} & {\rm RK} & {\rm Taylor}    \\[2ex]
5.85\, 10^{-10} & 6.30\,10^{-4} & 5.43\, 10^{-10} & 2.60\,10^{-10}  
\end{array}
$}

\bigskip

TABLE 6.- Errors with respect to the exact solution of the second iterations using Extended Picard (first column) and VIM (second column). The third and forth columns represent the errors due to the use of Runge-Kutta and Taylor of order ten, respectively.

\vskip1cm

Thus, in this particular case, the precision of Extended Picard  looks better than the precision of VIM and of the order of Runge-Kutta and Taylor of order ten.

\section{Some applications in chemical models}

Once we have introduced the  modification of the Picard method in order to obtain approximate solutions to first order non-linear systems of ordinary differential equations, it is time to use it in concrete applications with interest in Chemistry, as well as give arguments to show some of the advantages of this method. Two are the applications we want to discuss in this article. One is an analysis on the solutions of the equation of the glycolisis. The other one finds approximate solutions for the Brusselator, a system of equations that controls chemical reactions.  

\subsection{Glycolisis}

The glycolisis is a chemical reaction consisting in  the breakdown of glucose by enzymes, releasing energy and pyruvic acid and governs the glucose metabolism in humans and animals. A mathematical approach for a quantitative analysis of this reaction has been proposed in \cite{STR,SEL,DON}. See also \cite{BR,BR1,BB}. In two dimensions, a qualitative analysis yields to the following system of differential equations:

\begin{equation}\label{50}
y'= -y + f(y,z)\,, \qquad z'= b- f(x,y)\,, \qquad f(y,z)= az+ zy^2\,, \qquad a,b>0\,,
\end{equation}
where $y(t)$ is the concentration of adenosine diphosphate (ADP) and $z(t)$ is the concentration of fructuose-6-phosphate (F6P), $t$ being time. The derivative in \eqref{50} is taken with respect to time. Our objective is to obtain approximate solutions for this model of glycolisis reaction based in the method described above. Some qualitative properties of \eqref{50} are considered in \cite{STR}, in particular, the existence of a unique fixed point with coordinates $(y^*,z^*)=(b,b/(a+b^2))$, as well as the limit cycle and the Hopf bifurcation, which is described on the parameter space in terms of the constants.

In order to determine the latter, let us consider the linear approximation of \eqref{50} on a neighbourhood of the fixed point, which has the form $\mathbf y'=A\mathbf y$, where $\mathbf y=(y(t),z(t))^T$, where the $T$ means transpose and $A$ is the Jacobian matrix of the transformation in right hand side of \eqref{50}. The determinant of this matrix is $\det A= a+b^2>0$ and its trace, $\tau$, is 

\begin{equation}\label{51}
\tau= -\frac{b^4+(2a-1)b^2 +a(1+a)}{\det A}\,,
\end{equation}
which shows that the fixed point is unstable for $\tau>0$ and stable for $\tau<0$. In the parameter space with coordinates $(a,b)$, the fixed point is either stable or unstable depending on which region of this parameter lies. These two regions are separated by a curve of equation (this curve is determined by two functional expressions that we distinguish with signs plus or minus):

\begin{equation}\label{52}
b_\pm = \frac 1{\sqrt 2} \,\sqrt{1-2a \pm \sqrt{1-8a}}\,.
\end{equation} 

In the region determined by the values of $a$ and $b$ by $0<a<1/8$ and $b_+<b<b_+$, the fixed point is unstable. On the other hand, stable limit cycles do exists, which may be determined numerically. Outside that region, the fixed point is asymptotically stable.  The system shows a Hopf bifurcation as $b$ increases and $a<1/8$  and fixed. These results are similar to those discussed in \cite{STR}. 

This discussion on the cyclic limit shows that \eqref{50} is a stiff equation. This implies that, when integrating  the solution referred to one coordinate, say $y(t)$ or $z(t)$, we may observe certain time intervals in which they are some somehow unexpected variations on the solution, see Figure  6 (same for Brusselator in Figure 7). This makes it often difficult the integration process. With this method, we have managed to integrate the equation up to a good degree of accuracy. 

Coming back to \eqref{50}. As in \eqref{3}, we may write the iterative system, which in the present case takes the form, where we have omitted the dependence of the time variable for simplicity,

\begin{equation}\label{53}
y'_k= -y_k +a\, z_k + z_{k-1}\,y^2_{k-1}\,, \qquad z'_k= b-a\, z_k - z_{k-1}\, y^2_{k-1}\,.
\end{equation}

Let us compare \eqref{53} to that obtained with the standard Picard method:

\begin{equation}\label{54}
y'_k= -y_{k-1} +a\, z_{k-1} +z_{k-1}\,y^2_{k-1}\,,  \qquad z'_k= b- a\, z_{k-1} - z_{k-1}\,y^2_{k-1}\,.
\end{equation}

Let $[\alpha,\beta]$ be the interval of integration on the independent variable t, the time and assume we have divided into subintervals of equal length, as in the general case. For practical reasons, we shall usually choose the interval $[\alpha,\beta]$ equal to $[0,1]$, unless otherwise stated.

For a qualitative analysis, we need to specify the values of the parameters in \eqref{50}. Let us choose $a=0.4$ and $b=0.6$. For these values, system \eqref{50} has an asymptotically stable fixed point\footnote{In order to get a limit cycle, we need to choose smaller values of some of the parameters, say $a=0.04$ and $b=0.6$.}. The concentrations of ADP and F6P, given by $y(t)$ and $z(t)$, respectively, oscillate with $t$ up to reach the fixed point. 

To apply the method, we need to fix initial values that we fix as $y(\alpha)=1$ and $z(\alpha)=1$ and an initial seed that plays the role of $\mathbf y_0(t)=(y_0(t),z_0(t))$, which in our case is $y_0(t)\equiv z_0(t)=1$, for all values of $t$ runing out the interval $[0,1]$\footnote{This is an obvious choice for the interval $[\alpha,\beta]$}.  

It is noteworthy that just with two iterations, the extended Picard method gives a similar precision to the obtained using the traditional Runge-Kutta (RK) method, as proved by routine numerical test. 

However, should we enlarge the interval of integration, say by considering an interval of the form $[0,\beta]$ with $\beta>>1$, the bigger the distance to $0$ within the interval, the slower  the convergence of iterations to the exact solution, see Figure 1.  In order to solve this inconvenience, one may increase the number of iterations. This obvious remedy not only greatly increases the CPU times, but also produces spurious solutions due to a stockpile of numerical errors. This is why we need to combine the extended Picard with the segmentary integration, as was explained in Section 2.2. Let us consider the inhomogeneous terms $z_{k-1}\, y^2_{k-1}$ and $b- z_{k-1}\, y^2_{k-1}$, which appear in \eqref{51} and \eqref{52} and $-y_{k-1} +a\, z_{k-1} +z_{k-1}\,y^2_{k-1}$ and $b- a\, z_{k-1} - z_{k-1}\,y^2_{k-1}$ in \eqref{52}. We approximate these terms using minimal squares with polynomials of order one and three on each segment.  The combination of both methods greatly improves the precision, although it is certain that increases CPU times.

\begin{figure}
\centering
\includegraphics[width=0.5\textwidth]{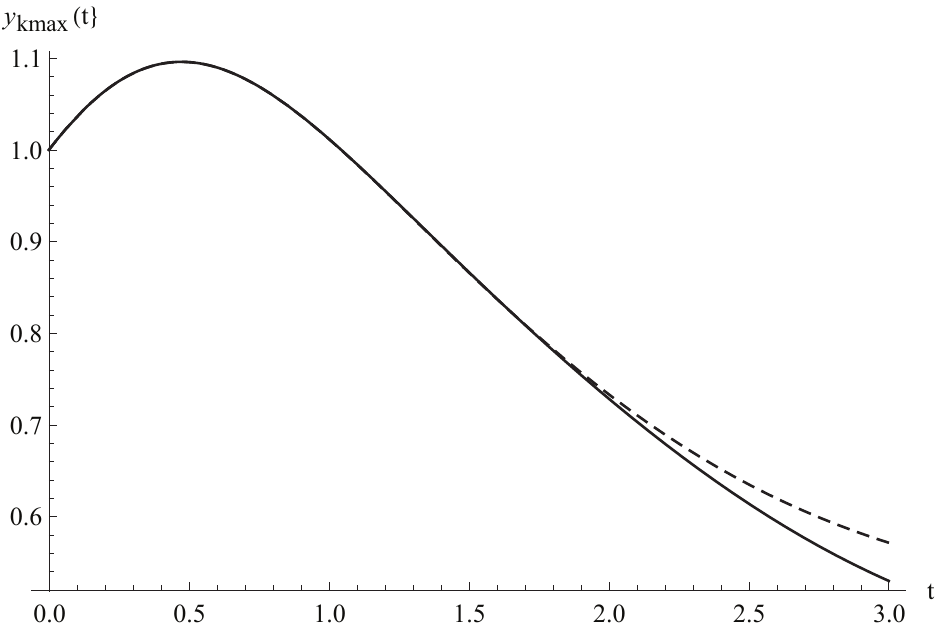}
\caption{\small Local solution to the Glycolisis equation for ADP concentration by iteration (dashed line) without making use of segmentary iteration, as compared with the Runge Kutta solution, when using $a=0.4$ and $b=0.6$ in \eqref{48}. We have made use of two iterations only, yet we have a reasonable fit for at least the half of the interval.  Note that both curves fit quite well up to values of $t$ of the order $t=1.5$ or bigger and, then, both curves split. This split is what the segmentary iteration avoids.
\label{Figure4}}
\end{figure}

In order to check the efficiency of the method, it is necessary to compute the error that the method has in relation with the exact solution. If this exact solution is not computable, then, we should obtain this error in comparison with some other method of reference \cite{NOR}. This is what we are going to do here, where the method of reference we choose is the RK of eighth order, for which its solution is here denoted as $y_{RK8}$. To evaluate the error, we use a standard formula such as

\begin{equation}\label{55}
e:= \frac1\beta \int_0^\beta ( y_{RK8}- y(t))^2\,dt\,,
\end{equation}
where $y(t)$ is the first component of $\mathbf y(t)=(y(t),z(t))$, the approximate solution we want to compare to $y_{RK8}$. This could be the approximate solution obtained by the Extended Picard, Standard Picard or some other like Taylor. We have evaluated the error obtained using these three methods in various numerical experiments. We here present some particular results that, nevertheless, show the general pattern of these experiments. 

In Figure 5, we give the solution for the ADP concentration using segmentary iteration. It is noteworthy that our solution exactly matches with the Runge-Kutta solution, at least on the time interval $[0,10]$. In Figure 6, we have chosen $a=0.04$ and $b=0.6$. Then, the fixed point becomes unstable. We have used segmentary iteration. On each interval, we have used three iterations and the approximation for the non-linear terms in the inhomogeneous part of the equation has been made using minimal squares with polynomials of degree three. We have used a much longer interval $[0,40]$. Note that the first part of the curve is similar to that shown in Figure 5.  

\begin{figure}
\centering
\includegraphics[width=0.5\textwidth]{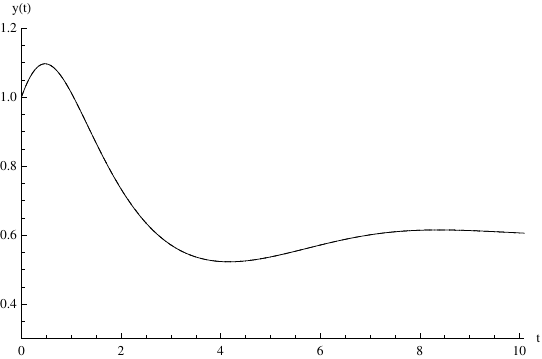}
\caption{\small Local solution to the Glycolisis equation for ADP concentration by  segmentary iteration up to $t=10$.  It completely matches  with the Runge Kutta solution. We have taken \eqref{21} $a=0.4$ and $b=0.6$. Note than in this case the fixed point is asymptotically stable.
\label{Figure5}}
\end{figure}

\begin{figure}
\centering
\includegraphics[width=0.5\textwidth]{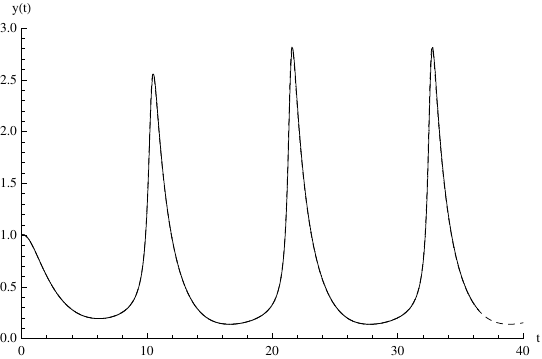}
\caption{\small Local solution to the Glycolisis equation for ADP concentration by  segmentary iteration up to $t=40$. We have taken \eqref{21} $a=0.04$ and $b=0.6$, so that the limit cycle exists. It completely matches  with the Runge Kutta solution, where the dashed line at the end denotes the extension of the RK solution. 
\label{Figure6}}
\end{figure}

The interval of integration chosen, $[0,40]$ is  wide enough to a proper test of Extended Picard. We divide this interval into subintervals of the same length $h=0.1$. To begin with, let us obtain the relative error to use the Taylor method as compared to the Runge-Kutta of eighth order. Here, $T_n$ means  Taylor method of order $n$. A simple table of errors gives:

\vskip1cm

\centerline{$
\begin{array}
[c]{cccc}
T_2 & T_3 & T_4  & T_5  \\[2ex]
2.5\,10^{-7} & 7.1\, 10^{-10} & 1.3\,10^{-12} & 9.6 \,10^{-15} 
\end{array}
$}

\medskip
TABLE 7.- Relative errors obtained with the use of the Taylor integration method of orders $n=2,3,4,5$ as compared with the Runge-Kutta method of order 8.

\vskip 1cm

It is also quite interesting to compare the precision obtained with the Extended Picard with the precision for Standard Picard. This is shown in the next two tables. The former refers to Extended Picard while the second one to Standard Picard. Let us go back to the comments after \eqref{53} and \eqref{54}, where we have stated that on each of the subintervals we approximate the inhomogeneous terms using minimal squares of order one and three. Numbers 1 and 3 on the top of Tables 8 and 9 refer to these orders. We consider up to five iterations on each case.

\vskip1cm

\centerline{$
\begin{array}
[c]{ccc}
{\rm iteration} & 1 & 3  \\[2ex]
2 & 3.6\, 10^{-11} & 3.8\,10^{-12}  \\
3 & 3.7\, 10^{-11}  & 6.2\, 10^{-15}  \\
4  & 3.7\, 10^{-11} & 6.0\, 10^{-15} \\
5 & 3.7\,10^{-11} & 6.0 \,10^{-15}
\end{array}
$}

\medskip
TABLE 8.- Relative errors obtained with the use of the Extended Picard method with iterations $n=2,3,4,5$ as compared with the Runge-Kutta method of order 8.

\vskip 1cm

\centerline{$
\begin{array}
[c]{ccc}
{\rm iteration} & 1 & 3  \\[2ex]
2 & 2.5\, 10^{-10} & 8.8\,10^{-10}  \\
3 & 5.2\, 10^{-10}  & 1.3\, 10^{-12}  \\
4  & 5.3\, 10^{-10} & 5.0\, 10^{-15} \\
5 & 5.2\,10^{-10} & 6.3 \,10^{-15}
\end{array}
$}

\medskip
TABLE 9.- Relative errors obtained with the use of the Standard Picard method with iterations $n=2,3,4,5$ as compared with the Runge-Kutta method of order 8.

\vskip 1cm

At first, let us compare Tables 8 and 7. We see that the Extended Picard with three iterations  gives the same degree of precision than a Taylor method of fifth order. In order to get the same precision for the Standard Picard, one needs one more iteration. It is in order to say that for the Runge-Kutta method the local error is estimated to be $o(h^n)$, where $n$ is the order of the method (8 in our case). 

As one may expect, for larger values of the subinterval width $h$, we have worst results in general. Nevertheless, there is some general behaviour that could be summarized as:

1.- Extended and Standard Picard up to five iterations give better results than Taylor up to order five.

2.- For values of $h=1$ or similar, Extended Picard gives better results than Standard Picard for 2 or 3 iterations. However, for more than 3 iterations Standard Picard may give better results. However, such a big subinterval width is never used for a simple reason: On the one step methods like ours the local error is of the order of $o(h^p)$, where $p$ is the method order. 

\subsection{Brusselator}

The Brusselator is a theoretical model for a kind of reaction called autocatalytic, which was developed by Prigogine and collaborators \cite{PL,GP,LE}. In two dimensions, it has the following form:

\begin{equation}\label{56}
y'(t) =1-(1+b)y(t) +f(y,z)\,,\qquad z'(t) = by-f(yz)\,,\qquad f(y,z)=ay^2z\,,
\end{equation}
with $a,b>0$ and also $y(t),z(t)>0$. 

This system has the following properties: It has a fixed point at $(y^*,z^*)=(1, b-a-1)$. The determinant of the Jacobian matrix evaluated at the fixed point is $1>0$ and its trace is $\tau=b-a-1$. The fixed point remains stable for $\tau<0$ and unstable for $\tau>0$. The $(a,b)$ plane is divided into two sectors separated by the curve $b=1+a$. On the region $\{a>0,b<1+a\}$ the fixed point is asymptotically stable and outside unstable. Outside the region of asymptotical stability of the fixed point we can still determine stable limit cycles using numerical experiments. 

Then, let us proceed to iterative integration. For any $k=1,2,\dots$, we have the equations

\begin{equation}\label{57}
y'_k =-(1+b) y_k +1 + a\, z_{k-1}\,y^2_{k-1}\,,\qquad z'_k = b\,y_k -a\, z_{k-1}\,y^2_{k-1}\,.
\end{equation}

If we choose the values $a=1$ and $b=2.5$ for the parameters, equation \eqref{56} shows a limit cycle on the plane $(y,z)$. Let us take initial values in the neighbourhood of the values determining this limit cycle, say $y_0=1.8$ and $z_0=1.2$. As time interval, we choose $[0,15]$ and for the width of each segment we take $h=0.1$. We approximate the inhomogeneous terms in \eqref{57}, $1 + a\, z_{k-1}\,y^2_{k-1}$ and $-a\, z_{k-1}\,y^2_{k-1}$, on each segment of width $h=0.1$ by minimal squares with polynomials of orders one aor three, which gives a polynomial particular solution. 

As done in previous cases, we estimate the relative error \eqref{55} with respect to a eighth order Runge-Kutta. Nevertheless, we first estimate the relative error with respect to Runge-Kutta of a traditional one step method as is the Taylor method. We consider here, second, third, forth and fifth order Taylor. The result is

\vskip1cm 

\centerline{$
\begin{array}
[c]{cccc}
T_2 & T_3 & T_4 & T_5  \\[2ex]
2.7\,10^{-5} & 4.6\, 10^{-7} & 9.4\,10^{-9} & 2.7\, 10^{-10} 
\end{array}
$}

\medskip
TABLE 10.- Relative errors, with respect to Runge-Kutta, obtained using the Taylor method with orders 2-5.

\vskip1cm

The relative errors obtained with the use of Extended Picard are shown on Table 11:

\vskip1cm

\vskip1cm

\centerline{$
\begin{array}
[c]{ccc}
{\rm iteration} & 1 & 3  \\[2ex]
2 & 3.5\, 10^{-4} & 1.1\,10^{-6}  \\
3 & 3.3\, 10^{-4}  & 1.0\, 10^{-7}  \\
4  & 3.3\, 10^{-4} & 8.8\, 10^{-8} \\
5 & 3.3\,10^{-4} & 8.8 \,10^{-8}
\end{array}
$}

\medskip
TABLE 11.- Errors of Extended Picard relative to Runge-Kutta. The number of iterations varies from 2 to 5 and the numbers 1 and 3 on top mean the degree of the polynomial used in the approximation by minimal squares.

\vskip1cm 

The repetition of the figures on Table 10 is only apparent and is due to our choice of only two digits. Deviations arise when we choose three or more digits. The figures obtained when using the Standard Picard are, however, similar. 

Nevertheless, with a simple trick we may improve the precision of Extended Picard for \eqref{56}. Let us use the new variable $w: = y+z$. Using this new variable $w$, we may transform system \eqref{56} into a second order differential equation of the form

\begin{equation}\label{58}
w''(t)+ a\,w(t)= F(w,w')\,,
\end{equation}
with

\begin{equation}\label{59}
F(w,w') = a+b -(1+b+a(3-2w))w' -a(w-3)w'^2 - a w'^3\,.
\end{equation}

Then, let us write \eqref{58} as a first order system by performing the change of variables given by $y:=w$, $z:=w'$. Needless to say that these $y$ and $z$ are different from the variables denoted by the same letters in \eqref{56}. Nevertheless, we keep this notation, as no confusion should arise. Then, \eqref{58} is transformed into

\begin{equation}\label{60}
y'=z\,,\qquad z'=-a\,y +F(y,z)\,,
\end{equation}
where $F(y,z)$ comes obviously from \eqref{59}. The values of the parameters $a=1$ and $b=2.5$ give a limit cycle, so that we choice these values and then integrate \eqref{60} using Extended Picard on the interval $[0,15]$ with $h=0.1$ as before.  Just as before, we evaluate the relative error using Extended Picard with respect to Runge-Kutta. Just as before, we first evaluate the relative error using second, third, forth and fifth order Taylor with respect to Runge Kutta. These errors are given in Table 12:

\vskip1cm

\centerline{$
\begin{array}
[c]{cccc}
T_2 & T_3 & T_4 & T_5  \\[2ex]
1,8\,10^{-5} & 5.1\, 10^{-7} & 2.5\,10^{-9} & 3.1\, 10^{-11} 
\end{array}
$}

\medskip
TABLE 12.- Relative errors, with respect to Runge-Kutta, obtained using the Taylor method with orders 2-5.

\bigskip

We see that these errors are quite similar to those obtained for  \eqref{56}. However, integration of \eqref{60} by Extended Picard just as done before for \eqref{56} gives the following results:

\vskip1cm

\centerline{$
\begin{array}
[c]{ccc}
{\rm iteration} & 1 & 3  \\[2ex]
2 & 33.1\, 10^{-5} & 3.8\,10^{-5}  \\
3 & 31.5\, 10^{-6}  & 1.9\, 10^{-7}  \\
4  & 6.5\, 10^{-7} & 6.0\, 10^{-10} \\
5 & 6.9\,10^{-7} & 6.3 \,10^{-13}
\end{array}
$}

\medskip
TABLE 13.- Values of error for \eqref{58} under the same conditions than those given on TABLE 11.

\vskip1cm

We see that the precision has neatly improved. This also shows that the precision of the method depends on the system under study. 

In Figure 7, we show the solution for the $y(t)$ component of the Brusselator, equation \eqref{27}. We have used a segmentary iteration, with two iterations and have approximated the non-linear inhomogeneous terms by minimal squares using polynomials of degree three. The considered time interval has been $[0,15]$.

\begin{figure}
\centering
\includegraphics[width=0.5\textwidth]{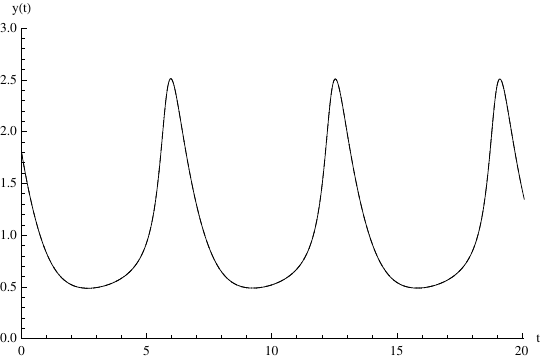}
\caption{\small Local solution corresponding to $y(t)$ for equation \eqref{54} where we use  segmentary iteration with two iterations and and minimal square approximation of order three. Here, the time interval is $[0,20]$. Observe the presence of peaks which are a clear sign of the stiffness of the equation. 
\label{Figure7}}
\end{figure}

\section{Concluding Remarks}

We analyze a method to study first order non-linear ordinary systems of equations that may appear in some mathematical models useful to describe chemical reactions. We call it the {\it Extended Picard method}. This method gives approximate analytic solutions, which are local, being given an initial or boundary condition. Depending on the particular equation under the study, these solutions may be directly find with the extension of Picard method without making use of the segmentation, which is a great advantage in relation with numerical methods which give discrete solutions. Approximate solutions are obtained by successive iterations, similarly to the standard Picard method, to which the procedure here introduced generalizes.  With respect to the precision of the method here display, we have compare it with the precision given by other simple one step methods such as Taylor.

 In some of the cases studied by ours, it seems necessary to perform a high number of iterations in order to have a good approximation. The higher the number of iterations the closer that the approximate solution is to the exact solution. It provides local solutions and this may have some disadvantages. In fact, for long integration intervals and large number of iterations, undesirable numerical  errors often occur, as illustrated in Figure 4. The cure is a combination of the iteration method with a segmentary integration. The whole interval is divided into smaller intervals, on which we apply the iteration method. Approximate solutions match  at the nodes, so that the approximate solution on the whole interval is continuous. For large integration intervals, non-linear terms in the inhomogeneous term of the equation are replaced by a minimal squares approximation using polynomials of degree one or three. Thus, we need much less iterations and avoid random errors. 

We have shown that the sequence of successive iterative solutions converge uniformly to the exact solution on compact sets of the real line. 

We have tested the precision of the method on three classic equations: Mathieu, quintic Duffing and Bratu.

We have used extended Picard  for two systems describing chemical reactions: Glycolisis and Brusselator equations. We have analyzed the behaviour of solutions and compared the relative precision of the method compared with Runge-Kutta with the precision of Taylor or Standard Picard, and we found it satisfactory. It is noteworthy that the precision  may depend on the specific system under study, as we have confirmed in the analysis of two different forms of the system describing the Brusselator. 

The procedure here studied may be extended to differential systems and equation of higher order. We shall attempt this research in a near future.

\section*{Acknowledgements}

The paper has been partially supported by the Q-CAYLE project, funded by the European Union-Next Generation UE/MICIU/Plan de Recuperacion, Transformacion y Resiliencia/Junta de Castilla y Leon (PRTRC17.11), and also by RED2022-134301-T, PID2020-113406GB-I00 and PID2023-148409NB-I00, financed by MICIU/AEI/10.13039/501100011033.

 \bigskip

{\bf Data Availability Statement}: All necessary data are in the paper.

\end{document}